\documentclass[reqno,11pt]{amsart} 
\usepackage{amsmath, amsfonts, amsthm, amssymb, verbatim,dsfont}
\usepackage{graphicx}
\usepackage[mathcal]{euscript} 
\usepackage{amstext} 
\usepackage{amssymb,mathrsfs}  
\theoremstyle{plain}
        \newtheorem{theorem}{Theorem}[section]
        \newtheorem{proposition}[theorem]{Proposition} 
        \newtheorem{lemma}[theorem]{Lemma}
        \newtheorem{corollary}[theorem]{Corollary} 
        \newtheorem{definition}[theorem]{Definition} 
\theoremstyle{definition}
        \newtheorem{remark}[theorem]{Remark}  
\theoremstyle{plain}
          
\numberwithin{equation}{section} 
\newcommand \be           {\begin{equation}}
\newcommand \ee            {\end{equation}}
  
\newcommand \Ecal       	{\mathcal{E}}  
  
\newcommand \RR           {\mathbb{R}}

\newcommand \del           \partial
\newcommand \eps        	\epsilon 
\newcommand \ubar   	{{\overline u}}

\newcommand \loc 		{{\mathrm{loc}}} 

\DeclareMathOperator	\sgn  {sgn}

\newcommand \la 		\langle
\newcommand \ra 		\rangle

\newcommand \nn       	{\mathbf n}

\newcommand \Ocal  	{\mathcal O}

\newcommand \TT		{\mathcal{T}}
\newcommand \Tcal	{\TT}
\newcommand{\dK}{\partial K}

\newcommand{\Lip}{\mathop{\mathrm{Lip}}}
\newcommand{\unk}{u^n_K}

\newcommand{\Hcal}{\mathcal{H}}

\newcommand{\Lorentz}{\mathbf{M}^{d+1}}
\newcommand{\Lorentzplus}{\mathbf{M}^{d+1}_+}
\DeclareMathOperator\divex {div}  
\newcommand{\dive}{\divex_{g}}
\DeclareMathOperator\diam {diam}

\newcommand{\Kn}{K^{n}}

\newcommand{\ukn}{{u_{K}^{n}}}
\newcommand{\ukne}{u_{K_{\ekz}}^{n}}
\newcommand{\uknn}{{u_{K}^{n+1}}}

\newcommand{\ekp}{e_{K}^{+}}
\newcommand{\ekm}{e_{K}^{-}}
\newcommand{\ekpm}{e_{K}^{\pm}}
\newcommand{\ekz}{e^{0}}
\newcommand{\ez}{\ekz}

\newcommand{\qq}{{\mathbf{q}_{K,\ekz}}}
\newcommand{\qqn}{{\mathbf{q}_{\Kn,\ekz}}}

\newcommand{\ukezm}{u_{K_{e^{0}}}^{-}}
\newcommand{\ukp}{u_{K}^{+}}
\newcommand{\ukm}{u_{K}^{-}}

\newcommand{\mukp}{\mu_{K}^{+}}
\newcommand{\mukm}{\mu_{K}^{-}}

\newcommand{\mukezp}{\overline\mu_{K,e^{0}}^{+}}

\newcommand{\Hkez}{H_{K,e^{0}}}
\newcommand{\sumezdk}{\sum_{e^{0}\in\del^{0}K}}
\newcommand{\QQ}{\mathbf{Q}_{K,\ekz}}
\newcommand{\QQn}{{\mathbf{Q}_{\Kn,\ekz}}}

\newcommand{\mukn}{ \mu_{K}^{n}}
\newcommand{\muknn}{ \mu_{K}^{n+1}}
\newcommand{\sumkez}{\substack{\Kn\in\Kcal^{n}\\ \ekz\in\del^{0}\Kn}}
\newcommand{\sumnN}{\sum_{n=0}^{N}}
\newcommand{\Kcal}{\mathcal{K}}
\newcommand{\Mbf}{\mathbf{M}} 

\newcommand{\wbf}{\mathbf{w}} 
\newcommand{\pkp}{ p^{+}_{K}}
\newcommand{\pkz}{ p^{0}_{K}}
\def\Xint#1{\mathchoice 
{\XXint\displaystyle\textstyle{#1}}%
{\XXint\textstyle\scriptstyle{#1}}%
{\XXint\scriptstyle\scriptscriptstyle{#1}}%
{\XXint\scriptscriptstyle\scriptscriptstyle{#1}}%
\!\int} 
\def\XXint#1#2#3{{\setbox0=\hbox{$#1{#2#3}{\int}$} 
\vcenter{\hbox{$#2#3$}}\kern-.5\wd0}} 
 
\def\dashint{\Xint\diagup}

\begin{document}

\title[Finite volume schemes on Lorentzian manifolds]
{Finite volume schemes on Lorentzian manifolds}
\author
%
[\textsc{P. Amorim, P.G. L{\tiny e}Floch,} \and \textsc{B. Okutmustur}]
{\textsc{Paulo Amorim, Philippe G. L{\tiny e}Floch,} \and \textsc{Baver Okutmustur}}
\address
{\textsc{P. Amorim, P.G. LeFloch,} \and \textsc{B. Okutmustur}\\
Laboratoire Jacques-Louis Lions \& Centre National de la Recherche Scientifique
\\
Universit\'e de Paris 6, 4 Place Jussieu, 75252 Paris, France.
\newline
E-mail : \tt Amorim@ann.jussieu.fr, LeFloch@ann.jussieu.fr, Okutmustur@ann.jussieu.fr.}
 
\date{\today}



\date{July 15, 2007} 

\medskip
\begin{abstract} 
We investigate the numerical approximation of (discontinuous) entropy solutions to 
nonlinear hyperbolic conservation laws posed on a Lorentzian manifold. 
Our main result establishes the convergence of monotone and first-order finite volume schemes 
for a large class of (space and time) triangulations. 
The proof relies on a discrete version of entropy inequalities and an entropy dissipation bound, 
which take into account the manifold geometry accurately
and generalize techniques and estimates that were known in the (flat) Euclidian setting, only.
The strong convergence of the scheme then is then a consequence of
the well-posed theory recently developed by Ben-Artzi and LeFloch for conservation laws on manifolds. 

\end{abstract}
\maketitle


\section{Introduction}
\label{IN}

In the present paper, we consider discontinuous solutions to nonlinear hyperbolic conservation laws 
posed on a globally hyperbolic Lorentzian manifold. Our main objective is to introduce
a class of first-order and monotone finite volume schemes and derive geometrically natural 
and nonlinear stability properties satisfied by these schemes. In turn, we will 
conclude that the proposed finite volume schemes converge (in a strong topology) 
toward the entropy solutions recently defined in Ben-Artzi and LeFloch \cite{BL}, who 
established the well-posedness theory for conservation laws posed on (Riemannian or Lorentzian) manifolds. 
Our proof can be regarded as a generalization to Lorentzian manifolds of 
the method introduced in Cockburn, Coquel and LeFloch \cite{CCL} for the (flat) Euclidean setting
and extended to Riemannian manifolds in Amorim, Ben-Artzi, and LeFloch \cite{ABL}.  

Major difficulties arise in working with a partial differential equation on a Lorent\-zian manifold, 
rather than on the customary (flat) Euclidian space.  A space and time triangulations must be 
introduced and the geometry of the manifold must be taken into account
accurately in the discretization.  
In addition, one cannot canonically choose a preferred time foliation in general, so that 
it is important for the discretization to be robust enough to allow for a large class 
of foliation and of space and time triangulations, that assume only limited regularity.
From the numerical analysis standpoint, it is challenging to design and analyze schemes 
that are indeed consistent with the geometry of the given Lorentzian manifold.  
 
Our assumption of global hyperbolicity of the Lorentzian background allows us to encompass 
physically realistic situation (arising, for instance, in the applications to general relativity). 
The class of schemes proposed in the present paper is quite general, and 
essentially assumes a monotonicity property on the flux functions only. 
Moreover, we over also a large class of space and time triangulations in which element can 
degenerate in certain directions of space or time.  

Several new difficulties appear when trying to generalize to Lorentzian manifolds 
the convergence results in \cite{ABL,CCL}. 
Our guide in deriving the necessary estimates was to ensure that all of our arguments 
are intrinsic in nature and (as the physics or geometry imposes it) 
do not rely explicitly on a choice of local coordinates. 

We show here that the proposed finite volume schemes can be expressed as a convex decomposition,  
and we derive a discrete version of entropy inequalities as well as 
sharp estimates on the entropy dissipation. Strong convergence towards an entropy solution
follows from DiPerna's uniqueness theorem \cite{DP}.

For further work on conservation laws on manifolds, we refer the reader to a pioneering paper by 
Panov \cite{P}, as well as to 
LeFloch and Okutmustur \cite{LO} and LeFloch, Neves, and Okutmustur \cite{LNO}.

An outline of this paper follows. In Section~\ref{CL} we briefly state some results 
about conservation laws posed on manifolds. In Section~\ref{FVM} we present
a formal derivation of the finite volume schemes and state our assumptions and 
main convergence result.  
Section \ref{PC} is devoted to establishing stability estimates which are of independent interest
but are also the basis to the convergence proof presented in Section \ref{POC}. 
Finally, in Section \ref{EX}, we present some examples of discretizations and particular schemes.


\section{Preliminaries} 
\label{CL}

\subsection{Conservation law on a Lorentzian manifold}
We will need some results on the well-posedness of hyperbolic conservation laws 
posed on Lorentzian manifolds., for which we refer to \cite{BL}.

Let $(\Lorentz,g)$ be a time-oriented, $(d+1)$-dimensional Lorentzian manifold. Here, 
$g$ is a metric with signature $(-, +, \ldots, +)$, and we denote by $T_p\Lorentz$  the tangent space 
at a point $p \in \Lorentz$.
Recall that tangent vectors $X$ (and more generally vector fields) 
on a Lorentzian manifold can be separated into time-like vectors
($g(X,X) < 0$), null vectors ($g(X,X) = 0$), and space-like vectors ($g(X,X) > 0$).
The manifold is time-oriented, i.e.~a consistent orientation can be chosen throughout the manifold,
so that  we can distinguish between past-oriented and future-oriented vectors. 
To the metric $g$ one associates its Levi-Cevita connection (covariant derivative) $\nabla$, 
which allows us to define the divergence $\dive$ operator, as in the Riemannian case, as the trace of the covariant derivative 
of a vector field. By duality arguments, 
the divergence of vector fields is known to depend only on the volume form associated with $g$. 
However, to assume a Lorentzian metric is most convenient to apply the terminology arising in general relativity
when dealing with the initial-value problem (see below).  

Recall that, by definition, a \emph{flux} on the manifold $\Lorentz$
is a vector field $f(\ubar, p) \in T_p\Lorentz$ depending on a real parameter $\ubar$.
The \emph{conservation law} on $(\Lorentz,g$) associated with a flux $f$ is
\be
\dive \big( f(u, p) \big) = 0, \qquad u: \Lorentz \to \RR.
\label{CL.1}
\ee
It is said to be \emph{geometry compatible} if $f$ satisfies the condition
\be
\dive f (\ubar, p) = 0, \quad \ubar \in \RR, \, p \in \Lorentz.
\label{CL.2}
\ee
Furthermore, $f$ is said to be a \emph{time-like flux} if
\be
g\big(\del_u f(\ubar, p), \del_u f(\ubar, p) \big) < 0, \qquad p \in \Lorentz, \, \ubar \in \RR.
\label{CL.3}
\ee

We are interested in the initial-value problem associated with \eqref{CL.1}. We fix a space-like hypersurface
$\Hcal_0 \subset \Lorentz$ and a measurable and bounded function $u_0$ defined on $\Hcal_0$.
Then, we search for a function $u=u(p) \in L^\infty(\Lorentz)$ satisfying \eqref{CL.1}  in the distributional
sense and such that the (weak) trace of $u$ on $\Hcal_0$ coincides with $u_0$:
\be
u_{|\Hcal_0} = u_0.
\label{CL.4}
\ee
It is natural to require that the vectors $\del_u f(\ubar, p)$, which determine the propagation of waves in
solutions of  \eqref{CL.1}, are time-like and future-oriented.

We assume that the manifold $\Lorentz$ is \emph{globally hyperbolic,} in the sense that
there exists a foliation of $\Lorentz$ by space-like, compact, oriented hypersurfaces $\Hcal_t$ ($t \in \RR$):
$$
\Lorentz = \bigcup_{t \in \RR} \Hcal_t.
$$
Any hypersurface $\Hcal_{t_0}$ is referred to as a \emph{Cauchy surface} in $\Lorentz$,
while the family $\Hcal_t$ ($t \in \RR$) is called an \emph{admissible foliation associated with} $\Hcal_{t_0}$.
The future of the given hypersurface will be denoted by
$$
\Lorentzplus := \bigcup_{t \geq 0} \Hcal_t.
$$
Finally we denote by $n^t$ the future-oriented, normal vector field to each $\Hcal_t$,
and by $g^t$ the induced metric.
Finally, along $\Hcal_t$, we denote by $X^t$ the normal component of a vector field $X$, thus
$X^t := g(X,n^t)$.

\begin{definition}
\label{CL-0}
 A flux $F=F(\ubar, p)$ is called a
\emph{convex entropy flux} associated with the conservation law \eqref{CL.1}
if there exists a convex function $U:\RR \to \RR$ such that
$$
F(\ubar, p) = \int_0^\ubar \del_u U(u') \,\del_u f(u', p) \, du', \qquad p \in \Lorentz, \, \ubar \in \RR.
$$
A measurable and bounded function $u=u(p)$ is called an \emph{entropy solution}
of conservation law \eqref{CL.1}--\eqref{CL.2}
if the following \emph{entropy inequality}
$$
\aligned 
& \int_{\Lorentzplus} g(F(u), \nabla_{g} \phi) \, dV_g + \int_{\Lorentzplus} (\dive F)(u) \, \phi \, dV_g 
\\
& + \int_{\Hcal_0} g_0(F(u_0), n_0) \, \phi_{\Hcal_0} \, dV_{g_0}  -\int_{\Lorentzplus} U'(u) (\dive f)(u)\, \phi\, dV_{g} \geq 0
\endaligned 
$$
holds for all convex entropy flux $F=F_p(\ubar)$ and all smooth functions $\phi \geq 0$ compactly supported in
$\Lorentzplus$.
\end{definition}

In particular, the requirements in the above definition imply the inequality
$$
\dive \big( F(u) \big) - (\dive F)(u)  + U'(u) (\dive f)(u) \leq 0
$$
in the distributional sense.


\subsection{Well-posed theory}

We will use the following result (see \cite{BL}): 

\begin{theorem}
\label{CL-1}
Consider a conservation law \eqref{CL.1}
posed on a globally hyperbolic Lorentzian manifold $\Lorentz$ with compact slices. 
Let $\Hcal_0$ be a Cauchy surface in $\Lorentz$, and $u_0: \Hcal_0 \to \RR$ be a function in 
$L^\infty(\Hcal_0)$. Then, the initial-value problem \eqref{CL.1}--\eqref{CL.4}
admits a unique entropy solution $u=u(p) \in L_\loc^\infty(\Lorentz_+)$.
Moreover, for every admissible foliation $\Hcal_t$, the trace $u_{|\Hcal_t} \in L^1(\Hcal_t)$ exists as a Lipschitz continuous
function of $t$. When the flux if geometry compatible, the functions
$$
\| F^t(u_{|\Hcal_t}) \|_{L^1(\Hcal_t)},
$$
are non-increasing in time, for any convex entropy flux $F$. Moreover, given any two entropy solutions $u,v$,
the function
$$
\| f^t(u_{|\Hcal_t}) - f^t(v_{|\Hcal_t}) \|_{L^1(\Hcal_t)}
$$
is non-increasing in time.
\end{theorem}

Throughout the rest of this paper, we therefore assume that a globally hyperbolic Lorentzian manifold is given, 
and we tackle the problem of the discretization 
of the initial value problem associated with the conservation law \eqref{CL.1} 
and a given initial condition where $u_0 \in L^\infty(\Hcal_{0})$.  
We \emph{do not} assume that the flux $f$ is geometry compatible,  
and on the map $\dive f$ we require the following \emph{growth condition}: there exist constants $C_1,C_2>0$ 
such that for all $(\ubar,p) \in \RR \times M$ 
\be
\label{FVM.2}
|(\dive f )(\ubar,p)| \le C_1 + C_2 \, |\ubar|. 
\ee

Two remarks are in order. 
First of all, the terminology here differs from the one in the Riemannian (and Euclidean) cases, 
where the conservative variable is singled out. The class of conservation laws on a Riemannian manifold is recovered 
by taking $\Lorentz= \RR \times \bar M$, where $\bar M$ is a Riemannian manifold 
and $f(\ubar, p ) = (\ubar, \bar f) \in \RR\times T_pM$, which leads to
$$
\dive \big( f(u) \big) = \del_{t} u + \dive \big( \bar f(u) \big).
$$

Second, in the Lorentzian case, no time-translation property is available in general, contrary to the
Riemannian case. Hence, no time-regularity is implied by the $L^1$ contraction property.
 

\section{Notation and main results}
\label{FVM} 

\subsection{Finite volume schemes}

Before we can state our main result we need to introduce several notations and motivate the formulation of the 
finite volume schemes under consideration.

We consider a (in both space and time) triangulation $\TT^{h} = \bigcup_{K\in\TT^{h}} K$, on $\Lorentz$ composed of bounded space-time elements $K$ satisfying the following assumptions:
\begin{itemize}
\item The boundary $\dK$ of an element $K$ is a piecewise smooth $d$-dimensional manifold without boundary, $\dK = \bigcup_{e\subset \dK} e$, where each $d$-dimensional face $e$ is a smooth manifold and is either everywhere time-like, null, or space-like. The outward unit normal to $e\in\del K$ is denoted by $\nn_{K,e}$.
\item Each element $K$ contains exactly two space-like faces, with disjoint interiors, denoted $\ekp$ and $\ekm$, such that the outward unit normals  to $K$, $\nn_{K,\ekp}$ and $\nn_{K,\ekm}$, are future- and past-oriented, respectively. They will be called the \emph{inflow} and the \emph{outflow faces}, respectively.
\item For each element $K$, the set of the \emph{lateral faces} $\del^{0}K :=\dK \setminus \{\ekp , \ekm\}$ is nonempty and time-like.
\item For every pair of distinct elements $K, K' \in \TT^{h}$, 
the set $K \cap K'$ is either a common face of $K, K'$ or else a submanifold with 
dimension at most $(d-1)$. 
\item $\Hcal_{0} \subset \bigcup_{K\in\TT^{h}} \dK$, where $\Hcal_{0}$ is the initial Cauchy hypersurface.
\item For every ${K\in\TT_{h}}$, $\diam \ekpm \le h$.
\end{itemize}
Given an element $K$, we denote by $K^{+}$ (resp. $K^{-}$) the unique element distinct from $K$ sharing the face $\ekp$ (resp. $\ekm$) with $K$, and for each $\ekz \in \del^{0}K$, we denote $K_{\ekz}$ the unique element sharing the face $\ekz$ with $K$. 

The most natural way of introducing the finite volume method is to view the discrete solution as defined on the space-like faces $e^{\pm}_{K}$ separating two elements. So, to a particular element $K$ we may associate two values, $\ukp$ and $\ukm$ associated to the unique outflow and inflow faces $\ekp,$ $\ekm$. Then, one may determine that the value $u_{K}$ of the discrete solution on the element $K$, is the solution $\ukp$ determined on the inflow face $\ekm$ (one could just as well say that $u_{K}$ is the solution $\ukp$ determined on the outflow face $\ekp$, or some average of the two, as long as one does this coherently throughout the manifold). 

Thus, for any element $K$, integrate the equation \eqref{CL.1}, apply the divergence theorem and decompose the boundary $\dK$ into its parts $\ekp, \ekm$, and $\del^{0}K$:
\be
\label{FVM.3}
\aligned
& - \int_{\ekp} g_{p}(f(u,p), \nn_{K,\ekp}(p)) dp - \int_{\ekm} g_{p}(f(u,p), \nn_{K,\ekm}(p)) dp
\\
&\hspace{43mm} + \sum_{e^{0}\in\del^{0}K}\int_{e^{0}} g_{p}(f(u,p), \nn_{K,e^{0}}(p)) dp = 0.
\endaligned
\ee
Note the minus signs on the first two terms. Indeed, for a Lorentzian manifold, the divergence theorem reads
$$
\int_{\Omega} \dive f dV_{\Omega} = \int_{\del\Omega} g(f , \tilde\nn) dV_{\del\Omega},
$$
with $\tilde \nn$ is the outward normal if it is space-like, and the inward normal if it is time-like.

Given an element $K$, we want to compute an approximation $\ukp$ of the average of $u(p)$ in the outflow face $\ekp$, being given the values of $\ukm$ on $\ekm$ and of $\ukezm$ for each $e^{0}\in\del^{0}K$.

The following notation will be useful. Let $f$ be a flux on the manifold $\Lorentz$, $K$ an element of the triangulation, and $e \subset \del K$. Define the function $\mu_{K,e}^{f} : \RR \to \RR$ by
\be
\label{FVM.3.0}
\mu_{K,e}^{f} (u) := \dashint_{e} g_{p} \big(f(u,p), \nn_{K,e}(p) \big) dp.
\ee
Also, if $w :\Lorentz\to\RR$ is a real function, we write
$$
\mu^{w}_{e} := \dashint_{e} w(p) dp.
$$
Using this notation, the  second term in \eqref{FVM.3} is approximated by
$$ 
\int_{\ekm} g_{p}\big(f(u,p), \nn_{K,\ekm}(p) \big) dp \simeq |\ekm| \mu_{K,\ekm}^{f} (\ukm),
$$
and the last term is approximated using
$$ 
\int_{e^{0}} g_{p}(f(u,p), \nn_{e^{0}}(p)) dp \simeq |e^{0}| \, \qq(\ukm, \ukezm),
$$
where to each element $K$, and each face $e^{0}\in\del^{0}K$ we associate a locally Lipschitz \emph{numerical flux function} $\qq(u,v ) :\RR^{2}\to \RR$ satisfying the assumptions given in \eqref{FVM.3.1}--\eqref{FVM.3.3} below.

Therefore, in view of the above approximation formulas
we may write, as a discrete approximation of \eqref{FVM.3},
\be
\label{FVM.6}
|\ekp| \mu_{K^{+},\ekp}^{f}(\ukp) := |\ekm| \mu_{K,\ekm}^{f}(\ukm) - \sum_{e^{0}\in \del^{0}K}|e^{0}| \, \qq(\ukm, \ukezm),
\ee
which is the \textbf{finite volume method} of interest and, equivalently
\be
\label{FVM.6.1}
\ukp := (\mu_{K^{+},\ekp}^{f})^{-1} \Big(\frac{ |\ekm|}{|\ekp|} \mu_{K,\ekm}^{f}(\ukm) - \sum_{e^{0}\in \del^{0}K}\frac{ |e^{0}|}{|\ekp|} \, \qq(\ukm, \ukezm) \Big).
\ee
The second formula which may be carried out numerically 
(using for instance a Newton algorithm) is justified by the following observation: 

\begin{lemma}
\label{FVM-0-1}
For any $K\in\TT^{h}$, the function $u \mapsto \mu_{K,\ekm}^{f}(u)$ is monotone increasing.
\end{lemma}

\begin{proof}
From \eqref{FVM.3.0} we deduce that 
$$
\del_{u}\mu_{K,\ekm}^{f}(u) = \dashint_{e} g_{p} \big(\del_{u}f(u,p), \nn_{K,\ekm}(p) \big) dp > 0,
$$
since $\del_{u} f(u,p)$ is future-oriented and $\nn_{k,\ekm}$ is past-oriented.
\end{proof}
If $\ekm \subset \Hcal_{0}$, the initial condition \eqref{CL.4} gives
\be
\label{FVM.6.2}
\ukm := \mu_{\ekm}^{u_{0}} = \dashint_{\ekm} u_{0}(p) dp 
\ee

Finally, we define the function $u^{h} : \Lorentz \to \RR$ by
\be
\label{FVM.6.3}
u^{h}(p) := \ukm, \qquad p\in K.
\ee
For $e\in\del K$ we introduce the notation
\be
\label{FVM.2.1}
f_{e}(u,p) := g_{p}(f(u,p), \nn_{K,e}(p)).
\ee


\subsection{Assumptions on the numerical flux-functions} 

We are now in a position to state our main result. First, we need to present some assumption on the flux function
and the triangulations. 

The numerical flux functions $\qq(u,v) : \RR^{2} \mapsto \RR$ in the equation \eqref{FVM.6} satisfy the following assumptions:
\begin{itemize}
\item{\emph{Consistency property :}
	\be 
	\label{FVM.3.1}
	\qq(u,u) = \dashint_{\ekz} f_{\ekz}(u,p) dp = \mu^{f}_{K,\ekz}(u).
	\ee
	}
\item{\emph{Conservation property :} 
	\be
	\label{FVM.3.2}
	\qq(u,v) = -\mathbf{q}_{K_{e^{0}},e^{0}}(v,u ), \qquad u,v \in \RR. 
	\ee
	}
\item{\emph{Monotonicity property :}
	\be
	\label{FVM.3.3}
	{\del_{ u}} \qq(u,v) \ge 0, \quad {\del_{ v}}\qq(u,v)  \le 0.
	\ee
	}
\end{itemize}

\begin{remark}
In fact, as we shall see in the examples in Section~\ref{EX}, it is natural, and even unavoidable, that the numerical flux-functions $\qq$ depend on the whole geometry of the element $K$ as well as of its neighbour, $K_{\ekz}$. However, to keep the notations as simple as possible, we will use the above notations throughout.
\end{remark}

For each element $K$, define the local time increment
$$
\tau_{K} = \frac{|K|}{|\ekp|}.
$$
We suppose that $\tau := \max\limits_{K} \tau_{K} \to 0$ when $h\to 0$, and that for every $K$,
\be
\label{FVM.A}
h^{2}/ \tau \to 0.
\ee
For stability, we impose the following \emph{CFL condition}
which should hold for all $ K\in \TT^{h},\,  \ekz\in\del_{0} K$: 
\be
\label{FVM.7}
\frac{|\del^{0}K|}{|\ekp|}  \sup_{u\in \RR} \big|\del_{u} \mu^{ f}_{K,\ekz}(u) \big| \le \big(\sup_{u\in \RR} \del_{u}(\mu^{ f}_{K^{+},\ekp})^{-1}(u) \big)^{-1} <\infty,
\ee


\subsection{Assumptions on the triangulations}

We now describe the time-like geometric structure of the triangulation $\Tcal^{h}$. We will introduce a 
(global and geometric in nature) 
admissibility condition (see \eqref{FVM.E} below) involving only the time evolution of the triangulation, which will be completely independent of its structure on space-like related faces. We stress that our method poses almost no restriction on the space-like structure of the discretization.

We will define the \emph{local Cartesian deviation} associated with a pair of elements $K$, $K^{-}$, which is a quantity depending only on the geometry of each such pair of space-time elements, and propose an admissibility condition involving a global bound on this quantity. The local Cartesian deviation measures the amount by which the time evolution of the triangulation deviates from a uniform, Cartesian evolution, and is defined as follows. Let $K\in \Tcal^{h}$, and suppose first that its space-like outflow face $\ekp$ satisfies the property that its center of mass, denoted by $ p^{+}_{K}$, lies on $\ekp$. Consider also the center of mass of $\del^{0}K$, denoted by $ p^{0}_{K}$. Then, define the bilinear form $\Ecal(K)$ on $(T_{ \pkp} \Mbf)^{2}$ by
$$
\Ecal(K) := \frac{1}{\tau_{K}} \wbf_{K} \otimes \nn_{K,\ekp},
$$
where the vector $\wbf \in T_{\pkp}\Mbf$ is the future-oriented unit tangent vector at $\pkp$ to the unique (if $h$ is small enough) geodesic connecting $\pkp$ to $\pkz$. So, if $X,Y$ are two vectors defined at the point $ p^{+}_{K}$, we have, by definition,
$$
\Ecal(K) (X,Y) = \frac{1}{\tau_{K}}\wbf_{K} \otimes \nn_{K,\ekp} (X,Y) = \frac{1}{\tau_{K}} g_{ p^{+}_{K}}(X,\wbf_{K})\, g_{ p^{+}_{K}} ( Y, \nn_{K,\ekp}),
$$
and we define the \textbf{local Cartesian deviation} associated with $K$, $K^{-}$, as the form
$$
|K| \Ecal(K) - |K^{-}| \Ecal(K^{-}).
$$

This tensor-like operator measures the rate of change of the quantity 
$$|K| \Ecal(K)(X,Y)$$
with respect to the time-like direction defined locally by the faces $\ekpm$. Our admissibility criterion then states that this rate of change must not blow up faster than $h^{-1}$ after summation over $K\in\Tcal^{h}$, as well as some flatness conditions on the space-like faces $\ekpm$.

\begin{definition}
\label{FVM-0-2}
We say that $\Tcal^{h}$ is an \textbf{admissible triangulation} if the following conditions are satisfied: 
\enumerate
\item For  every smooth vector field $\Phi$ with compact support and every family of smooth vector fields $\Psi_{K}$, $K\in\Tcal^{h}$, the Cartesian deviation satisfies
\be
\label{FVM.E}
\Big| \sum_{\Kn\in\Tcal^{h}} \big( |K| \Ecal(K) - |K^{-}| \Ecal(K^{-}) \big) (\Phi, \Psi_{K}) \Big| \lesssim \frac{\eta(h)}{h}\, \| \Phi \|_{L^{\infty}} \sup_{K} \| \Psi_{K}\|_{L^{\infty}},
\ee
for some function $\eta(h)$ such that
$$
\eta(h) \to 0.
$$
\item For every $K\in\Tcal^{h}$, 
\be
\label{FVM.F}
d_{g}(\pkp , \ekp) \lesssim |\ekp| \tau,
\ee
where $d_{g}$ denotes the distance function associated to the metric $g$.
\item For all smooth vector fields $X$ defined on $\ekp$,
\be
\label{FVM.C}
\| g_{p}(X(p), \nn_{K,\ekp}(p)) \|_{C^{2}(\ekp)} \lesssim \| X(p) \|_{C^{2}(\ekp)}.
\ee
\endenumerate
\end{definition}

In the case where $\pkp \not\in \ekp$, the condition \eqref{FVM.F} means that as $h\to0$, the center of mass $\pkp$ of $\ekp$ approaches $\ekp$, even when we blow it up by a factor of order $|\ekp|$, as one would expect to find if $\ekp$ becomes flat in the limit. Under this assumption, we may uniquely extend the normal vector field $\nn_{K,\ekp}$ by geodesic transport to a small neighbourhood $V$ of $\ekp$, in such a way that $\pkp \in V$. Thus, in formula \eqref{FVM.E}, we can meaningfully speak of the vector $\nn_{K,\ekp}(\pkp) \in T_{\pkp}\Mbf$. The condition \eqref{FVM.C} is intended to rule out oscillations on the normal vector field due to the geometry of $\ekp$.

The assumption \eqref{FVM.E} amounts to a global geometric condition on the triangulation's Cartesian deviation. As we shall see in Section~\ref{EX}, this allows us to cover a large class of interesting, implementable space-time triangulations.


\subsection{Main result of convergence}

Finally, we are in a position to state: 

\begin{theorem}[\textbf{Convergence of the finite volume method}]
\label{FVM-1}
Let $u^{h}$ be the sequence of functions generated by the finite volume method \eqref{FVM.6}--\eqref{FVM.6.3} on an admissible triangulation, with initial data $u_{0}\in  L^\infty(\Hcal_{0})$, and with numerical flux-functions satisfying the conditions \eqref{FVM.3.1}--\eqref{FVM.3.3}, and the CFL condition \eqref{FVM.7}. Then, for every 
slice $\Hcal_t$, the sequence $u^h$ is bounded in $L^\infty\big(\bigcup_{s \in [0,t]} \Hcal_s\big)$, and converges almost everywhere when $h\to 0$ towards the unique entropy solution $u \in 
L^\infty(\Lorentz_+)$ of the Cauchy problem \eqref{CL.1}, \eqref{CL.4}.
\end{theorem}

In Section~\ref{PC} below, we will derive the key estimates required for a proof of Theorem \ref{FVM-1}
which will be the content of Section~\ref{POC}. We follow here the strategy originally proposed by Cockburn, Coquel and LeFloch 
\cite{CCL}. New estimates are required here to take into account the geometric effects and, especially, 
during the time evolution in the scheme. We will start with local (both in time and in space) entropy estimates, 
and next deduce a global-in-space entropy inequality. We will also establish the $L^{\infty}$ stability of the scheme
and, finally, the global (space-time) entropy inequality required for the convergence proof.

\begin{remark}
Our formulas do coincide with the formulas already known in the Riemannian and Euclidean cases. In these cases, the function $\mu^{f}_{K,\ekm}(u)$ is always the identity function. Therefore, the finite volume scheme reduces to the usual formulation found in \cite{CCL,ABL}. Also, our expression for the time increment $\tau$ and the CFL condition \eqref{FVM.7} give the usual formulas when particularized to the Euclidean or Riemannian setting.
\end{remark}

\begin{remark}
More generally, and in the interest of practical implementation, one may replace the right-hand side of the equations \eqref{FVM.3.0} and \eqref{FVM.6.2} with more realistic averages. For instance, one could  take an average of $g(f(u,p), \nn_{K,e})$ over $N$ spatial points $p_{j}$ given from some partition $e^{j}$ of $e$,
$$
\mu_{K,e}^{f}(u) = \frac{1}{|e|}\sum_{j=1}^{N}|e^{j}| g(f(u,p_{j}), \nn_{K,e}(p_{j})).
$$
In view of these remarks, we see that the finite volume method may be given more generally by the algorithm which consists of fixing an averaging operator $\mu_{K,\ekm}^{f}$, and using the equation \eqref{FVM.6.1} to iterate the method,
with initial data given by
$$
\ukm := \mu_{\ekm}^{ u_{0}}.
$$
However, any such average is just an approximation of the integral expression used in \eqref{FVM.3.0}. 
This approximation can be chosen to be of very high order on the parameter $h$, by choosing appropriate quadrature formulas. Therefore, for the sake of clarity, we will present the proofs with $\mu^{f}_{K,e}$ defined by \eqref{FVM.3.0} and omit the (straightforward) treatment of the error term issuing from this approximation.
\end{remark}


\section{Discrete entropy estimates}
\label{PC}

\subsection{Local entropy dissipation and entropy inequalities}

First of all, let us introduce some notations which will simplify the statement of the results as well as the proofs. Define
$$
\mukp(u) := \mu^{f}_{K^{+},\ekp}(u) = -\mu^{f}_{K,\ekp}(u)
$$
$$
\mukm(u) := \mu^{f}_{K,\ekm}.
$$
With this notation, the finite volume method \eqref{FVM.6.1} reads
\be
\label{FVM.7.0.1}
|\ekp| \mukp(\ukp) = |\ekm| \mukm(\ukm) - \sumezdk |\ekz| \qq(\ukm, \ukezm).
\ee
As in \cite{CCL} and \cite{ABL}, we rely on a convex decomposition of $\mukp(\ukp)$, which allows us to control the entropy dissipation. 

Define $\tilde \mu_{K,e^{0}}^{+}$ by the identity
$$
\tilde \mu_{K,e^{0}}^{+} :=  \mukp(\ukm) - \frac{|\del^{0}K|}{|\ekp|} \big( \qq(\ukm, \ukezm) - \qq(\ukm,\ukm)\big),
$$
and define 
\be
\label{FVM.7.2}
\mukezp :=  \tilde \mu_{K,e^{0}}^{+} - \frac{1}{|\ekp|} \int_{K} \dive f (\ukm,p) dp
\ee
Then, one has the following convex decomposition of $\mukp(\ukp)$, whose proof is immediate from \eqref{FVM.7.0.1}.
\be
\label{FVM.7.3}
\mukp(\ukp) = \frac{1}{|\del^{0}K|} \sumezdk |e^{0}| \mukezp.
\ee
 
\begin{lemma}
\label{PC-1}
Let $(U(u),F(u,p))$ be a convex entropy pair (cf. Definition \ref{CL-0}). For each $K$ and each $e=\ekm, \ekp$, let $V_{K,e} : \RR \to \RR$ be the convex function defined by 
\be
\label{FVM.7.3.1}
V_{K,e} (\mu) := \mu^{F}_{K,e}\big( (\mu^{f}_{K,e})^{-1}(\mu)\big), \quad \mu\in\RR.
\ee
Then there exists a family of numerical flux-functions $\QQ(u,v) : \RR^{2}\to\RR$ satisfying the following conditions:
\begin{itemize}
\item $\QQ$ is consistent with the entropy flux $F$:
$$
\QQ(u,u ) = \mu^{F}_{K,\ekz}(u),  \qquad K\in\Tcal^{h}, \ekz\in \del^{0}K, u\in \RR.
$$
\item Conservation property:
$$
\QQ(u,v ) = -\mathbf{Q}_{K_{\ekz},\ekz}(v,u), \qquad u,v \in \RR.
$$
\item Discrete entropy inequality:
\be
\label{FVM.7.4}
\aligned
& V_{K^{+},\ekp} \big( \tilde\mu_{K,e^{0}}^{+} \big)  -  V_{K^{+},\ekp} \big(  \mukp(\ukm) \big) 
\\
&\qquad+ \frac{|\del^{0}K|}{|\ekp|} \big( \QQ(\ukm,\ukezm) - \QQ(\ukm,\ukm) \big) \le 0.
\endaligned 
\ee
\end{itemize}
\end{lemma}

\

At this juncture, we conclude from the inequality \eqref{FVM.7.4} that 
\be
\label{FVM.7.5}
\aligned
&  V_{K^{+},\ekp} ( \mukezp )  -  V_{K^{+},\ekp} \big(  \mukp(\ukm) \big) 
\\
&\qquad+ \frac{|\del^{0}K|}{|\ekp|} \big( \QQ(\ukm,\ukezm ) - \QQ(\ukm,\ukm) \big) \le  R_{K,\ekz}^{+},
\endaligned 
\ee
where $R_{K,\ekz}^{+}$ is given by
\be
\label{FVM.7.6}
R_{K,\ekz}^{+} :=  V_{K^{+},\ekp} ( \mukezp ) -  V_{K^{+},\ekp}( \tilde \mu_{K,e^{0}}^{+}).
\ee

\begin{proof}
First of all, note that using \eqref{FVM.7.3.1} we may write the inequality \eqref{FVM.7.4} equivalently as
\be
\label{FVM.7.7}
\aligned
& \mu^{F}_{K^{+},\ekp}\big( (\mu^{f}_{K^{+},\ekp})^{-1} (\tilde \mu_{K,e^{0}}^{+}) \big)  - \mu^{F}_{K^{+},\ekp}(  \ukm ) 
\\
&\qquad+ \frac{|\del^{0}K|}{|\ekp|} \big( \QQ(\ukm,\ukezm) - \QQ(\ukm,\ukm) \big) \le 0.
\endaligned 
\ee
Indeed, we have for instance
$$
\aligned
  V_{K^{+},\ekp} \big(  \mukp(\ukm) \big) &= \mu^{F}_{K^{+},\ekp}\big( (\mu^{f}_{K^{+},\ekp})^{-1}(\mukp(\ukm))\big) 
=  \mu^{F}_{K^{+},\ekp}(\ukm) 
\endaligned
$$

We begin by introducing the following operator. For $u,v \in \RR$, $\ekz\in \del^{0}K,$ let 
$$ 
 H_{K,e^{0}}(u,v) :=  \mukp(u) - \frac{|\del^{0}K|}{|\ekp|} \big( \qq(u,v) - \qq(u,u)\big) 
$$
We claim that $\Hkez$ satisfies the following properties:
\be
\label{FVM.9}
\frac{\del}{\del u} \Hkez(u,v)\ge 0,
\qquad 
\frac{\del}{\del v}\Hkez(u,v) \ge 0,
\ee 
\be
\label{FVM.11}
 \Hkez(u,u) =  \mukp(u).
\ee
The first and last properties are immediate. The second is a consequence of the CFL condition \eqref{FVM.7} and the monotonicity of the method. Indeed, from the definition of $H_{K,e^{0}}(u,v)$
 we may perform exactly the same calculation as in the proof of Lemma~\ref{PC-2} to prove that $\Hkez(u,v)$ is a convex combination of $\mukp(u)$ and $\mukp(v)$, which in turn are increasing functions. This establishes 
 the first inequality in \eqref{FVM.9}. 

We now turn to the proof of the entropy inequality \eqref{FVM.7.7}. Suppose first that 
\eqref{FVM.7.7} is already established for the Kruzkov's family of entropies $\overline U(u,\lambda)=| u - \lambda |$, $\overline F(u,\lambda,p) = \sgn(u - \lambda)( f(u,p) - f(\lambda,p)),$ $\lambda\in\RR$. In this case, the Kruzkov numerical entropy flux-functions are given by
$$
\overline{\mathbf{Q}}_{K,\ekz}(u,v,\lambda ) := \qq(u\vee\lambda, v\vee \lambda) - \qq(u \wedge \lambda, v \wedge \lambda ).
$$
It is easy to check that $\QQ$ satisfies the first two conditions of the lemma.

In fact, it is enough to prove the inequality \eqref{FVM.7.7} for 
Kruzkov's entropies only. Indeed, if $U$ is a smooth function which is linear at infinity, we have
(formally)
$$
\aligned
\frac{1}{2}\int_{\RR}\overline U(u,\lambda) U''(\lambda) d\lambda = \frac{1}{2} \int_{\RR}\overline U''(u,\lambda) U(\lambda) d\lambda = \frac{1}{2}\langle \delta_{\lambda=u}, U(\lambda) \rangle = U(u),
\endaligned
$$
modulo an additive constant. Similarly, if $(U,F)$ is a convex entropy pair, we obtain
$$
\frac{1}{2}\int_{\RR}\overline  F(u,\lambda,p) U''(\lambda)d \lambda = F(u,p).
$$
Since we shall prove an $L^{\infty}$ bound for our approximate solutions, we may suppose that the $u$ above varies in a bounded set $B\subset\RR$. Thus, we may apply the same reasoning with any function which is not linear at infinity, by changing it into a linear function outside $B$. This shows that we can obtain the inequality \eqref{FVM.7.7} for any convex entropy pair $(U,F)$ by first proving it in the special case of Kruzkov's entropies, multiplying by $U''(\lambda)/2$, and integrating. In that case, the numerical flux will be given by
$$
\QQ(u,v) = \frac{1}{2} \int_{\RR}\overline{\mathbf{Q}}_{K,\ekz}(u,v,\lambda) U''(\lambda) d\lambda.
$$
Again, this numerical flux satisfies the first two assumptions of the lemma, since they are inherited from the corresponding properties for the Kruzkov numerical flux $\QQ(u,v,\lambda)$.

Therefore, we now proceed to prove the inequality \eqref{FVM.7.7} for Kruzkov's family of entropies. This is done in two steps. First, we will show that
\be
\label{FVM.11.1}
\aligned
\mu^{\overline F}_{K^{+},\ekp}(\ukm,\lambda) &- \frac{|\del^{0}K|}{|\ekp|} \big( \overline{\mathbf{Q}}_{K,\ekz}(\ukm,\ukezm,\lambda) - \overline{\mathbf{Q}}_{K,\ekz}(\ukm,\ukm,\lambda) \big) 
\\
&=  H(\ukm\vee \lambda, \ukezm\vee \lambda) - H(\ukm\wedge \lambda, \ukezm\wedge \lambda).
\endaligned
\ee
Second, we will see that for any $u,v,\lambda\in\RR$, we have
\be
\label{FVM.11.2}
\aligned
H(u\vee\lambda, v\vee\lambda) -  H(u\wedge\lambda, v\wedge\lambda) \ge
\mu^{\overline F}_{K^{+},\ekp}\big( (\mukp)^{-1} (H(u,v)) , \lambda\big).
\endaligned
\ee
For ease of notation, we omit $K, \ekz$ from the expression of $H$. The identity \eqref{FVM.11.1} and the inequality \eqref{FVM.11.2} (with $u=\ukm, v=\ukezm$) combined give \eqref{FVM.7.7}, for Kruzkov's entropies.

To prove \eqref{FVM.11.1}, simply observe that
$$
\aligned
\mu^{\overline F}_{K^{+},\ekp}(\ukm,\lambda) &=  \sgn(\ukm - \lambda) \big( \mu^{f}_{K^{+},\ekp}(\ukm) - \mu^{f}_{K^{+},\ekp}(\lambda) \big)
\\
& =  \sgn(\mukp(\ukm) - \mukp(\lambda)) \big( \mukp(\ukm) - \mukp(\lambda) \big)
\\
& = \big|\mukp(\ukm) - \mukp(\lambda) \big|
\\
& =  \big( \mukp(\ukm) \vee \mukp(\lambda) - \mukp(\ukm) \wedge \mukp(\lambda) \big)
\\
&=  \big(\mukp(\ukm \vee \lambda) - \mukp(\ukm \wedge \lambda) \big).
\endaligned
$$
Here, we have repeatedly used that $\mukp$ is a monotone increasing function. The identity \eqref{FVM.11.1} now follows from the expressions of the Kruzkov numerical entropy flux, $\overline{\mathbf{Q}}_{K,\ekz}$, and of $H$.

Consider now the inequality \eqref{FVM.11.2}. We have
$$
\aligned
H(u\vee\lambda, v\vee\lambda) & - H(u\wedge\lambda, v\wedge\lambda)
\\
& \ge \big( H(u,v)\vee H(\lambda,\lambda) \big) -  \big( H(u,v)\wedge H(\lambda,\lambda) \big).
\endaligned
$$
This is a consequence of the fact that if $\varphi$ is an increasing function, then $\varphi(u\vee\lambda) = \varphi(u\vee\lambda) \vee \varphi(u\vee\lambda) \ge \varphi(u) \vee \varphi(\lambda)$, and \eqref{FVM.9}. Thus,
we have  
\begin{align*}
& H(u\vee\lambda, v\vee\lambda) -  H(u\wedge\lambda, v\wedge\lambda) 
\\
&\ge  \big| H(u,v) - H(\lambda,\lambda) \big|
\\
& = \big|  H(u,v) - \mukp(\lambda) \big| 
\\
& =  \sgn\big(  H(u,v) - \mukp(\lambda) \big) \big( H(u,v) -\mukp(\lambda) \big)
\\
& =  \sgn\big( (\mukp)^{-1}\big(  H(u,v)\big) - \lambda \big) 
        \big( \mukp\big( (\mukp)^{-1} \big( H(u,v) \big)\big) -\mukp(\lambda) \big)
\\
& = \mu^{\overline F}_{K^{+},\ekp}\big( (\mukp)^{-1}\big( H(u,v)\big), \lambda \big).
\end{align*} 
This establishes \eqref{FVM.11.2}. We now choose $u=\ukm, v = \ukezm$ in \eqref{FVM.11.2}, observe that $H_{K,\ez}(\ukm, \ukezm) = \tilde\mu_{K,\ez}^{+}$, and combine with \eqref{FVM.11.1} to obtain inequality \eqref{FVM.7.7} for Kruzkov's entropies. As described above, \eqref{FVM.7.7} will hold for all convex entropy pairs $(U,F)$. 

To conclude the proof, observe that the functions $V_{K,e}$ in \eqref{FVM.7.4} are indeed convex. 
This follows from twice differentiating the representation formula 
\be
\label{FVM.100}
V_{K,e}(\mu) = \int^{\mu} U'\big((\mu^{f}_{K,e})^{-1}(u) \big) du,
\ee
and using the convexity of $U$ and the monotonicity of $(\mu^{f}_{K,e})^{-1}$, for $e=\ekp,\ekm$. This completes the proof of Lemma~\eqref{PC-1}.
\end{proof}


\subsection{Entropy dissipation estimate and $L^\infty$ estimate}

We now discuss the time evolution of the triangulation. As we have said, the initial hypersurface $\Hcal_{0}$ is composed of inflow faces $\ekm$. We then define the hypersurfaces $\Hcal_{n}$, for $n>0$, by 
$$
\Hcal_{n} := \bigcup_{\ekm \subset \Hcal_{n-1}} \ekp,
$$
and set
$$
\Kcal^{n} := \big\{ K : \ekm \subset \Hcal_{n-1} ,\quad \ekp \subset \Hcal_{n} \big\}.
$$
An element of $\Kcal^{n}$ is denoted by $K^{n}$ when it is important to stress the time level $n$. Otherwise, when writing local in time estimates, we omit the time level $n$. 

Next, we introduce the following notations. For $K^{n}\in \Kcal^{n}$, write
$$
\mu^{n}_{K} := \mu^{-}_{K^{n}} = \mu^{f}_{K^{n},e_{K^{n}}^{-}}, \quad u^{n}_{K} := u^{-}_{K^{n}}
$$
so that for instance, $\mu_{K^{n}}^{+}(u_{K^{n}}^{+}) = \mu^{n+1}_{K}(u^{n+1}_{K})$.
Accordingly, we define 
\be
\label{FVM.11.5}
V^{n}_{K}(\mu) := V_{K^{n},e_{K^{n}}^{-}}(\mu),
\ee
where $V_{K,e}$ are the time-like entropy flux defined in Lemma \ref{PC-1}.

\begin{lemma}
\label{PC-2}
The finite volume approximations satisfy the $L^{\infty}$ bound
\be
\label{FVM.12}
\max_{\Kn\in\Kcal^{n}} |\ukn| \le \big( \max_{K^{0}\in \mathcal{K}_{0}} |u_{K}^{0}|  + C_{1}t_{n}\big) e^{C_{2}t_{n}}
\ee
for some constants $C_{1,2}\ge 0$, where 
\be
\label{FVM.12.1}
t_{n} := \sum_{j=0}^{n} \tau_{j} = \sum_{j=0}^{n} \max_{K^{j}\in\mathcal{K}_{j}} \frac{|K^{j}|}{|e_{K^{j}}^{+}|} < \infty.
\ee
\end{lemma}

\begin{proof}
First of all, observe that from the consistency condition \eqref{FVM.3.1}, the definition of $\mu^{f}_{K,e}$ in \eqref{FVM.3.0} and the divergence theorem, we have for any $u\in\RR$,
$$
\aligned
\int_{\Kn} \dive f(u, p) dp &= \int_{\del\Kn} g_{p} \big(f(u, p), \tilde\nn(p) \big) dp
\\
& = |e_{K^{n}}^{+}| \muknn (u) -  |e_{K^{n}}^{-}| \mukn (u) + \sum_{\ekz\in\Kn} |\ekz| \qqn(u,u)
\endaligned
$$
(recall that $\tilde \nn$ is the interior unit normal if it is time-like, and the exterior unit normal if it is space-like). Moreover, with the current notations, the finite volume scheme \eqref{FVM.6} reads
$$
|e_{K^{n}}^{+}| \muknn(\uknn) = |e_{K^{n}}^{-}| \mukn(\ukn) - \sum_{\ekz\in\Kn} |\ekz| \qqn(\ukn, u_{K_{\ekz}}^{n}).
$$
Combining these two identities gives
\be
\aligned
\label{FVM.13}
\muknn(\uknn) &=  \muknn(\ukn) -  \frac{1}{|e_{K^{n}}^{+}|}\int_{\Kn} \dive f(\ukn, p) dp
\\
&\quad - \sum_{\ekz\in\Kn} \frac{ |\ekz|}{|e_{K^{n}}^{+}|} \big(  \qqn(\ukn, u_{K_{\ekz}}^{n}) - \qqn(\ukn,\ukn) \big).
\endaligned
\ee
Next, we rewrite the right-hand side as follows: 
\be
\label{FVM.14}
\aligned
\muknn(\uknn) &=  ( 1 -  \sum_{\ekz\in \del^{0}\Kn}\alpha_{\Kn,\ekz}) \muknn(\ukn)
+ \sum_{\ekz\in \del^{0}\Kn}\alpha_{\Kn,\ekz} \muknn(\ukne) 
\\
&\quad- \frac{1}{|e_{K^{n}}^{+}|}\int_{\Kn} \dive f(\ukn, p) dp,
\endaligned
\ee
where 
$$
\alpha_{\Kn,\ekz} := \frac{ |\ekz|}{|e_{K^{n}}^{+}|} \frac{\qqn(\ukn, \ukne) - \qqn(\ukn, \ukn)}{\muknn(\ukn) - \muknn(\ukne)}.
$$
This gives a convex combination of $\muknn(\ukn)$ and $\muknn(\ukne)$. Indeed, on one hand
we have $$
\sum_{\ekz\in\Kn} \alpha_{\Kn,\ekz} \ge 0
$$
due to the monotonicity condition \eqref{FVM.3.3} and Lemma \ref{FVM-0-1}. On the other hand, the CFL condition \eqref{FVM.7} gives
$$
\aligned
\sum_{\ekz\in\Kn}\alpha_{\Kn,\ekz}
& <  \Big| \frac{\ukn - \ukne}{\muknn(\ukn) - \muknn(\ukne)} \Big|  ( \Lip (\muknn)^{-1})^{-1}
\\
&\le \Lip (\muknn)^{-1} / \Lip (\muknn)^{-1} = 1.
\endaligned
$$
Thus, we find
$$
\aligned
&\muknn(\uknn)  \ge \min\big(\muknn(\ukn), \min_{\ekz\in\Kn} \muknn(\ukne)\big)  - \frac{1}{|e_{K^{n}}^{+}|}\int_{\Kn} \dive f(u, p) dp
\\
&\muknn(\uknn) \le \max\big(\muknn(\ukn), \max_{\ekz\in\Kn} \muknn(\ukne)\big)- \frac{1}{|e_{K^{n}}^{+}|}\int_{\Kn} \dive f(u, p) dp.
\endaligned
$$
Applying the monotone increasing function $(\muknn)^{-1}$, we find
$$
\aligned
&\uknn  \ge \min\big(\ukn, \min_{\ekz\in\Kn} \ukne \big)  + \frac{\Lip (\muknn)^{-1}}{|e_{K^{n}}^{+}|}\int_{\Kn} |\dive f(\ukn, p)| dp
\\
&\uknn \le \max\big(\ukn, \max_{\ekz\in\Kn} \ukne\big) + \frac{\Lip (\muknn)^{-1}}{|e_{K^{n}}^{+}|}\int_{\Kn} |\dive f(\ukn, p)| dp,
\endaligned
$$
which in turn gives
$$
|\uknn| \le  \max_{\Kn\in\Kcal^{n}}|\ukn| +  \max_{\Kn\in\Kcal^{n}}\frac{\Lip (\muknn)^{-1}}{|e_{K^{n}}^{+}|} \int_{\Kn} |\dive f(\ukn, p)| dp.
$$
By induction we obtain
$$
|\uknn| \le  \max_{K^{0}\in\mathcal{K}_{0}}|u_{K}^{0}| + \sum_{j=0}^{n}\max_{K^{j}\in\mathcal{K}^{j}}\frac{\Lip (\mu_{K}^{j+1})^{-1}}{|e_{K^{j}}^{+}|} \int_{K^{j}} |\dive f(u_{K}^{j}, p)| dp.
$$
Now use the growth condition \eqref{FVM.2} on the last term,
$$
\aligned
& \sum_{j=0}^{n}\max_{K^{j}\in\mathcal{K}^{j}}\frac{\Lip (\mu_{K}^{j+1})^{-1}}{|e_{K^{j}}^{+}|}\int_{K^{j}} |\dive f(u_{K}^{j}, p)| dp 
\\
& \le \sum_{j=0}^{n}\max_{K^{j}\in\mathcal{K}^{j}}\frac{\Lip (\mu_{K}^{j+1})^{-1}}{|e_{K^{j}}^{+}|} |K^{j}| \big( C_{1} + C_{2} |u_{K}^{j}| \big)
\\
&  \le  \big( C_{1} t_{n} + C_{2} \sum_{j=0}^{n}\tau^{j}\max_{j} |u_{K}^{j}| \big).
\endaligned
$$
Here, the constants $C_{1,2}$ may change from line to line. We have also used the fact that
$$
\max_{K^{j}\in\mathcal{K}_{j}}\Lip (\mu_{K}^{j+1})^{-1} \le C,
$$
which is an easy consequence of the hypotheses on the flux $f$. 
The result now follows from a discrete version of Gronwall inequality 
(see \cite[Lemma 6.1]{ABL}). This completes the proof of Lemma \ref{PC-2}.
\end{proof}

Recall that if $V$ is a convex function, then its \emph{modulus of convexity} on a set $S$ is defined by
$$
\beta := \inf \big\{ V''(w) : w \in S \big\}.
$$

\begin{proposition} Let $V^{n}_{K}$ be defined by \eqref{FVM.7.3.1}, \eqref{FVM.11.5}, and let $\beta_{K}^{n}$ be the modulus of convexity of $V_{K}^{n}$.
Then, we have
\label{PC-3}
\be
\label{FVM.20}
\aligned
&\sum_{\Kn\in\Kcal^{n}} |e_{K^{n}}^{+}| V^{n+1}_K(\muknn(\uknn)) 
\\
&\qquad +  \sum_{\Kn, \ekz}\frac{\beta_{K}^{n+1}}{2} \frac{|\ekz| |e_{K^{n}}^{+}|}{|\del^{0}\Kn|}  \big| \overline\mu_{K,\ekz}^{n+1} - \muknn(\uknn) \big| ^{2} 
\\
&\le \sum_{\Kn\in\Kcal^{n}} |e_{K^{n}}^{-}| V_{K}^{n}(\mukn(\ukn)) 
\\
&\qquad + \sum_{\Kn\in\Kcal^{n}} \int_{\Kn} \dive F(\ukn,p) dp  +\sum_{\substack{\Kn\in\Kcal^{n}\\ \ekz\in\del^{0}\Kn}} \frac{|\ekz| |e_{K^{n}}^{+}|}{|\del^{0}\Kn|} R_{K,\ekz}^{n+1} 
\endaligned
\ee
\end{proposition}

\begin{proof}
Consider the discrete entropy inequality \eqref{FVM.7.5}. Multiplying by $\frac{|\ekz| |e_{K^{n}}^{+}|}{ |\del^{0}\Kn|}$ and summing in $\Kn\in\Kcal^{n}$, $\ekz\in\del^{0} \Kn$ gives
\be
\label{FVM.21}
\aligned
&\sum_{\sumkez} \frac{|\ekz| |e_{K^{n}}^{+}|}{ |\del^{0}\Kn|} V^{n+1}_K ( \overline\mu_{K,\ekz}^{n+1} )  -  \sum_{\Kn\in\Kcal^{n}}|e_{K^{n}}^{+}| V^{n+1}_K \big(  \muknn(\ukn) \big) 
\\
& \quad + \sum_{\sumkez}|\ekz| \big( \QQn(\ukn,\ukne ) - \QQn(\ukn,\ukn) \big)
\\
& \le  \sum_{\sumkez}\frac{|\ekz| |e_{K^{n}}^{+}|}{ |\del^{0}\Kn|}R_{K,\ekz}^{n+1}.
\endaligned 
\ee
Next, observe that the conservation property \eqref{FVM.3.2} gives
\be
\label{FVM.22}
\sum_{\sumkez}|\ekz| \QQn(\ukn,\ukne ) = 0.
\ee
Now, if $V$ is a convex function, and if $v =  \sum_{j}\alpha_{j} v_{j}$ is a convex combination of $v_{j}$, then an elementary result on convex functions gives
$$
V(v) + \frac{\beta}{2} \sum_{j}\alpha_{j} |v_{j}- v|^{2} \le \sum_{j} \alpha_{j} V(v_{j}).
$$
Now, apply this result with the convex combination \eqref{FVM.7.3} and with the convex function $V^{n+1}_{K}$, multiply by $\frac{|\ekz| |e_{K^{n}}^{+}|}{ |\del^{0}\Kn|}$ and sum in $\Kn\in\Kcal^{n}$, $\ekz\in\del^{0} \Kn$. Combine the result with \eqref{FVM.21}, \eqref{FVM.22} to obtain
\be
\label{FVM.23}
\aligned
&\sum_{\Kn\in\Kcal^{n}} |e_{K^{n}}^{+}| V^{n+1}_K (\muknn(\uknn) )  -  \sum_{\Kn\in\Kcal^{n}}|e_{K^{n}}^{+}| V^{n+1}_K \big(  \muknn(\ukn) \big) 
\\
& \qquad + \sum_{\Kn, \ekz}\frac{\beta_{K}^{n+1}}{2} \frac{|\ekz| |e_{K^{n}}^{+}|}{|\del^{0}\Kn|}  \big| \overline\mu_{K,\ekz}^{n+1} - \muknn(\uknn) \big| ^{2} 
\\
&\qquad - \sum_{\sumkez}|\ekz| \QQn(\ukn,\ukn) 
\le  \sum_{\sumkez}\frac{|\ekz| |e_{K^{n}}^{+}|}{ |\del^{0}\Kn|}R_{K,\ekz}^{n+1}.
\endaligned 
\ee
Finally, using the identity
\be
\label{FVM.24}
\aligned
& \int_{\Kn} \dive F(u, p) dp = \int_{\del\Kn} g_{p} \big(F(u, p), \tilde\nn(p) \big) dp
\\
& = |e_{K^{n}}^{+}| V^{n+1}_{K}( \muknn (u)) -  |e_{K^{n}}^{-}|V^{n}_{K}( \mukn (u) )+ \sum_{\ekz\in\Kn} |\ekz| \QQn(u,u)
\endaligned
\ee
(with $u=\ukn$) yields the desired result. This completes the proof of Proposition \ref{PC-3}.
\end{proof}

\begin{corollary}
\label{PC-4}
Suppose that for each $K\in\TT^{h}$, $e = e_{K}^{\pm}$, the function $V_{K,e}$ is strictly convex, and that, moreover, we have 
\be
\label{FVM.25}
\beta_{K}^{n} \ge \beta > 0,
\ee
uniformly in $K$ and $n$. Then one has the following global estimate for the entropy dissipation,
\be
\label{FVM.26}
\aligned
\sum_{n=0}^{N}\sum_{\sumkez} \frac{|\ekz| |e_{K^{n}}^{+}|}{ |\del^{0}\Kn|} \big| \overline\mu_{K,\ekz}^{n+1} - \muknn(\uknn) \big|^{2} = \Ocal (t_{N}),
\endaligned
\ee
where $t_{N}$ is defined in \eqref{FVM.12.1}.
\end{corollary}

\proof
Summing the inequality \eqref{FVM.20} for $n=0,\dots, N$, we observe that the first terms on each side of the inequality cancel, leaving only the terms with $n=0$ and $n=N$. Moreover, using the growth condition \eqref{FVM.2} on the divergence term gives
\be
\label{FVM.27}
\aligned
& \sumnN\sum_{\Kn, \ekz}\frac{\beta_{K}^{n+1}}{2} \frac{|\ekz| |e_{K^{n}}^{+}|}{|\del^{0}\Kn|}  \big| \overline\mu_{K,\ekz}^{n+1} - \muknn(\uknn) \big| ^{2} 
\\
&\le  \sum_{K^{0}\in\mathcal{K}_{0}} |e_{K^{0}}^{-}| | V^{0}_K(\mu^{0}_{K}(u^{0}_{K})) |
+ \sum_{K^{N+1}\in\mathcal{K}_{N+1}} |e_{K^{N}}^{+}| |V^{N+1}_K(\mu^{N+1}_{K}(u^{N+1}_{K}))|
\\
&\qquad +\sumnN\sum_{\Kn\in\Kcal^{n}} |\Kn| \big( C_{1}+ C_{2} |u_{K}^{j}| \big)  +\sumnN\sum_{\substack{\Kn\in\Kcal^{n}\\ \ekz\in\del^{0}\Kn}} \frac{|\ekz| |e_{K^{n}}^{+}|}{|\del^{0}\Kn|} R_{K,\ekz}^{n+1}.
\endaligned
\ee
The last term is estimated using \eqref{FVM.7.2}, \eqref{FVM.7.6}, and the growth condition \eqref{FVM.2}, yielding
$$
\aligned
\sumnN\sum_{\substack{\Kn\in\Kcal^{n}\\ \ekz\in\del^{0}\Kn}} \frac{|\ekz| |e_{K^{n}}^{+}|}{|\del^{0}\Kn|} R_{K,\ekz}^{n+1} 
&\le \sumnN \sum_{\Kn\in\Kcal^{n}}\Lip V^{n+1}_{K} \int_{\Kn} |\dive f(\ukn,p)| dp
\\
& \le \sumnN \sum_{\Kn\in\Kcal^{n}} |\Kn| \big( C_{1}+ C_{2} |u_{K}^{j}| \big).
\endaligned
$$
Here, we have used that $\Lip V^{n+1}_{K}$ is uniformly bounded, which is an easy consequence of the corresponding bounds for the flux $f$. The result now follows from \eqref{FVM.25} and the $L^{\infty}$ estimate in Lemma \ref{PC-2}, which allows us to bound uniformly all the terms on the right-hand side of \eqref{FVM.27}. Note however that this bound depends, of course, on the entropy $U$. This completes the proof of Corollary \ref{PC-4}.
\endproof


\subsection{Global entropy inequality in space and time}

In this paragraph, we deduce a global entropy inequality from the local entropy inequality \eqref{FVM.7.5}. This is nothing but a discrete version of the entropy inequality used to define a weak entropy solution.  
Given a test-function $\phi$ defined on $\Lorentz$ we introduce its averages
\begin{align}
\nonumber 
&\phi_{\ekz}^{n} := \dashint_{\ekz} \phi(p) dp,
\\
\nonumber 
&\phi_{\del^{0}K}^{n} := \sum_{\ekz\in\del^{0}\Kn} \frac{|\ekz|}{|\del^{0}\Kn|} \phi_{\ekz}^{n} = \dashint_{\del^{0}\Kn} \phi(p) dp. 
\end{align}
The following lemma is easily deduced from the corresponding result in 
the Euclidean space, by 
relying on a system of local coordinates. This result will be useful when analyzing the approximation.

\begin{lemma}
\label{PC-4-1}
Let $G : \Mbf \to \RR $ be a smooth function, and let $e$ be a submanifold of $\Mbf$ such that $\mathrm{diam}(e) \le h$. Then there exists a point $p_{e}$ (not necessarily in $e$) such that
$$
\Big| \dashint_{e} G(p) dV_{e} - G(p_{e}) \Big| \lesssim h^{2}  \| G\|_{C^{2}(e)}.
$$
\end{lemma}

We are now ready to prove the global discrete entropy inequality.

\begin{proposition}
\label{PC-5}
Let $(U,F)$ be a convex entropy pair, and let $\phi$ be a non-negative test-function. 
Then, the function $u^h$ given by \eqref{FVM.6.3} satisfies the global entropy inequality
\be
\label{FVM.31}
\begin{aligned}
& -\sum_{n=0}^{\infty}  \sum_{\Kn\in\Kcal^{n}} \int_{\Kn}\dive \big(F(\ukn,p) \phi(p)\big) dp 
\\
& \le - \sum_{n=0}^{\infty} \sum_{\sumkez} \frac{|\ekz|}{|\del^{0}\Kn|} |e_{K^{n}}^{+}| \phi_{\ekz}^{n} \big( V^{n+1}_{K}(\tilde\mu^{n+1}_{K,\ekz}) - V^{n+1}_{K}(\overline\mu^{n+1}_{K,\ekz}) \big)
\\
& \quad +\sum_{n=0}^{\infty} \sum_{\sumkez} \frac{|\ekz|}{|\del^{0}\Kn|} |e_{K^{n}}^{+}| (\phi_{\del^{0}K}^{n} - \phi_{\ekz}^{n}) V^{n+1}(\overline\mu^{n+1}_{K,\ekz})
\\
& \quad +\sum_{n=0}^{\infty} \sum_{\sumkez} \int_{\ekz} (\phi_{\ekz}^{n} - \phi(p)) 
F_{\ekz}(\ukn,p) dp
\\
& \quad -\sum_{n=0}^{\infty} \sum_{\Kn\in\Kcal^{n}}
         \Big( \int_{e_{K^{n}}^{-}} (\phi_{\del^{0}K}^{n-1} - \phi(p)) \, g( F(\ukn,p) , \nn_{K^{n},\ekm}) dp 
\\
& \hskip4.cm + \int_{e_{K^{n}}^{+}} (\phi_{\del^{0}K}^{n} - \phi(p)) \,g( F(\ukn,p) , \nn_{K^{n},\ekp})dp \Big)
\\
& \quad + \sum_{K\in\Kcal_{0}} \int_{e_{K}^{-}}\phi_{\del^{0}K}^{0} g_{p}(F(u_{K}^{0},p) , \nn_{K,\ekm}) dp.
\end{aligned}
\ee
\end{proposition}

\proof
From the local entropy inequalities \eqref{FVM.7.5}, we obtain
\be
\label{FVM.36}
\begin{aligned} 
& \sum_{\sumkez} \frac{|\ekz| |e_{K^{n}}^{+}|}{ |\del^{0}\Kn|} \phi_{\ekz}^{n} \Big( V^{n+1}_K ( \overline\mu_{K,\ekz}^{n+1} ) 
- V^{n+1}_K \big(  \muknn(\ukn) \big) \Big)
\\
& \quad + \sum_{\sumkez}|\ekz| \phi_{\ekz}^{n}\big( \QQn(\ukn,\ukne ) - \QQn(\ukn,\ukn) \big)
\\ 
& \le \sum_{\sumkez}\frac{|\ekz| |e_{K^{n}}^{+}|}{ |\del^{0}\Kn|} \phi_{\ekz}^{n} R_{K,\ekz}^{n+1}.
\end{aligned}
\ee 
Now, 
from the conservation property \eqref{FVM.3.2} we have
$$ 
\sum_{\sumkez}|\ekz| \phi_{\ekz}^{n} \QQn(\ukn,\ukne ) = 0
$$
Also, from the consistency property \eqref{FVM.3.1}, we find
$$ 
\aligned
& \sum_{\sumkez}|\ekz| \QQn(\ukn,\ukn ) =  \sum_{\sumkez} \phi_{\ekz}^{n} \int_{\ekz} F_{\ekz}(\ukn, p) dp 
\\
& =  \sum_{\sumkez} \int_{\ekz} \phi(p) F_{\ekz}(\ukn, p) dp 
    +  \sum_{\sumkez} \int_{\ekz} (\phi_{\ekz}^{n} - \phi(p)) F_{\ekz}(\ukn, p) dp.
\endaligned
$$
Next, we have
$$
\aligned
& \sum_{\sumkez} \frac{|\ekz| |e_{K^{n}}^{+}|}{ |\del^{0}\Kn|} \phi_{\ekz}^{n} V^{n+1}_K ( \overline\mu_{K,\ekz}^{n+1} ) 
\\ 
& = \sum_{\sumkez} \frac{|\ekz| |e_{K^{n}}^{+}|}{ |\del^{0}\Kn|} \phi_{\del^{0}K}^{n} V^{n+1}_K ( \overline\mu_{K,\ekz}^{n+1} ) 
  +\sum_{\sumkez} \frac{|\ekz| |e_{K^{n}}^{+}|}{ |\del^{0}\Kn|} (\phi_{\ekz}^{n} - \phi_{\del^{0}K}^{n} ) V^{n+1}_K ( \overline\mu_{K,\ekz}^{n+1} )
\\
&\ge \sum_{\Kn\in\Kcal^{n}} |e_{K^{n}}^{+}| \phi_{\del^{0}K}^{n} V^{n+1}_K ( \muknn(\uknn) ) 
      +\sum_{\sumkez} \frac{|\ekz| |e_{K^{n}}^{+}|}{ |\del^{0}\Kn|} (\phi_{\ekz}^{n} - \phi_{\del^{0}K}^{n} ) V^{n+1}_K ( \overline\mu_{K,\ekz}^{n+1} )
\endaligned
$$
Here, we have used the fact that for a convex function $V$ and a convex combination $v = \sum_{j} \alpha_{j} v_{j}$ one has
$$
V(v) \le \sum_{j} \alpha_{j} V(v_{j}),
$$
with the convex function $V^{n+1}_{K}$ and the convex combination \eqref{FVM.7.3}.
Also, 
\begin{multline*}
\sum_{\sumkez}\frac{|\ekz| |e_{K^{n}}^{+}|}{ |\del^{0}\Kn|} \phi_{\ekz}^{n} V^{n+1}_K \big(  \muknn(\ukn) \big) 
=   \sum_{K^{n}\in\Kcal^{n}}|e_{K^{n}}^{+}| \phi_{\del^{0}K}^{n} V^{n+1}_K \big(  \muknn(\ukn) \big).
\end{multline*}
Therefore, the inequality \eqref{FVM.36} becomes 
\be
\label{FVM.39}
\aligned 
& \sum_{\Kn\in\Kcal^{n}} \phi_{\del^{0}K}^{n} |e_{K^{n}}^{+}| \Big(
 V^{n+1}_{K}\big(\muknn(\uknn)) - V^{n+1}_{K}\big(\muknn(\ukn)) \Big)
\\
& \quad -\sum_{\sumkez} \int_{\ekz} \phi(p) F_{\ekz}(\ukn, p) dp
\\
& \le   \sum_{\sumkez}\frac{|\ekz| |e_{K^{n}}^{+}|}{ |\del^{0}\Kn|} \phi_{\ekz}^{n} R_{K,\ekz}^{n+1} - \sum_{\sumkez} \frac{|\ekz| |e_{K^{n}}^{+}|}{ |\del^{0}\Kn|} (\phi_{\ekz}^{n} - \phi_{\del^{0}K}^{n} ) V^{n+1}_K ( \overline\mu_{K,\ekz}^{n+1} )
\\
& \quad + \sum_{\sumkez} \int_{\ekz} (\phi_{\ekz}^{n} - \phi(p)) F_{\ekz}(\ukn, p) \, dp
\\
& =: A^{h} + B^{h} + C^{h}
\endaligned
\ee
The first term in \eqref{FVM.39} can be written as
$$ 
\aligned 
& \sum_{\Kn\in\Kcal^{n}} \phi_{\del^{0}K}^{n} |e_{K^{n}}^{+}| \Big(
 V^{n+1}_{K}\big(\muknn(\uknn)) - V^{n+1}_{K}\big(\muknn(\ukn)) \Big)
\\
& = \sum_{\Kn\in\Kcal^{n}} \int_{e_{K^{n}}^{+}} \phi(p)\, g\big( F(\uknn, p) - F(\ukn, p), \tilde\nn_{K^{n},e^{+}_{K}}(p) \big) \, dp
\\
& \quad + \sum_{\Kn\in\Kcal^{n}} \int_{ e_{K^{n}}^{+}} (\phi_{\del^{0}K}^{n} - \phi(p)) \,g\big( F(\uknn, p) - F(\ukn, p), \tilde\nn_{K^{n},e^{+}_{K}}(p) \big) dp.
\endaligned 
$$ 
Combining this result with the identity
\be
\nonumber
\aligned 
& \int_{K} \dive \big( F(u, p) \phi(p)\big) dp = \int_{\del K} \phi(p)\, g( F(u,p), \tilde\nn_{\del K}) dp 
\\
& =\int_{e^{+}_{K}} \phi(p)\, g( F(u,p), \tilde\nn_{ K,e^{+}_{K}}) dp + \int_{e^{-}_{K}} \phi(p)\, g( F(u,p), \tilde\nn_{ K,e^{-}_{K}}) dp
\\
& \quad + \sum_{\ekz\in\del^{0}K} \int_{\ekz} \phi(p) F_{\ekz}(u, p) dp
\endaligned
\ee
(with $u = \ukn$) and in view of \eqref{FVM.39} we see that 
\be
\nonumber 
\aligned 
& -\int_{\Kn} \dive \big( F(\ukn, p) \phi(p)\big) dp 
\\
& \le A^{h} + B^{h} + C^{h} 
\\
& \quad -\sum_{\Kn\in\Kcal^{n}} \Big(
 \int_{e_{K^{n}}^{+}} \phi(p)\, g\big( F(\uknn, p), \tilde\nn_{K^{n},e^{+}_{K}}(p) \big) dp + \int_{e^{-}_{K}} \phi(p)\, g( F(u,p), \tilde\nn_{ K,e^{-}_{K}}) dp \Big)
\\
& \quad - \sum_{\Kn\in\Kcal^{n}} \int_{ e_{K^{n}}^{+}} (\phi_{\del^{0}K}^{n} - \phi(p)) \,g\big( F(\uknn, p) - F(\ukn, p), \tilde\nn_{K^{n},e^{+}_{K}}(p) \big) \, dp.
\endaligned 
\ee
The inequality \eqref{FVM.31} is now obtained by summation in $n$. First, the (summed) terms $A,B,C$ give the first three terms on the right-hand side of \eqref{FVM.31}, and by discrete integration by parts in the last terms of the above inequality, we find the remaining terms. This completes the proof of Proposition~\ref{PC-5}.
\endproof


\section{Proof of convergence}
\label{POC}

This section contains a proof of the convergence of the finite volume method, 
and is based on the framework of measure-valued solutions to conservation laws,
introduced by DiPerna \cite{DP} and extended to manifolds by Ben-Artzi and LeFloch \cite{BL}. 
The basic strategy will be to rely on the discrete entropy inequality \eqref{FVM.31} as well as on 
the entropy dissipation estimate \eqref{FVM.26}, in order
to check that any Young measure associated with the approximate solution is a measure-valued solution to the Cauchy problem 
under consideration. In turn, by the uniqueness result for measure-valued solutions it follows that, in fact, 
this solution is the unique weak entropy solution of the problem under consideration.

Since the sequence $u^h$ is uniformly bounded in 
$L^\infty( \mathbf{M})$, we can associate a subsequence and a
Young measure $\nu : \mathbf M \to \mathrm{Prob}(\RR)$,
which is a family of probability measures in $\RR$ parametrized by $p \in \mathbf M$. 
The Young measure allows us to determine all weak-$*$ limits of composite functions $a(u^h)$, 
for arbitrary real continuous functions $a$, according to the following property :
\be
\label{POC.1}
a(u^h) \mathrel{\mathop{\rightharpoonup}\limits^{*}} \langle \nu, a \rangle
\quad \text{ as } h \to 0
\ee
where we use the notation 
$$
\langle \nu, a \rangle := \int_\RR a(\lambda) \, d\nu(\lambda).
$$

In view of the above property, the passage to the limit in the 
left-hand side of \eqref{FVM.31} is (almost) immediate. The uniqueness theorem \cite{DP, BL} tells us that once we know that $\nu$ is a 
measure-valued solution to the conservation law, we can prove that the support
of each probability measure $\nu_{p}$ actually reduces to a single value
$u(p)$, if the same is true on $\Hcal_{0}$, that is, $\nu_{p}$ is the Dirac measure $\delta_{u(p)}$. 
It is then standard to deduce that the convergence in \eqref{POC.1} is actually strong, 
and that, in particular, $u^h$ converges strongly to $u$ which in turn is the unique 
entropy solution of the Cauchy problem under consideration. 

\begin{lemma}
\label{POC-1}
Let $\nu_{p}$ be the Young measure associated with the sequence $u^{h}$. Then, for every convex entropy pair $(U,F)$ and every non-negative test-function $\phi$ defined on $\Mbf$ with compact support, we have
\begin{multline}
\label{POC.5}
-\int_{\Mbf} \langle \nu_{p},\dive F(\cdot, p)\rangle \phi(p) + g\big( \langle \nu_{p}, F(\cdot,p)\rangle, \nabla \phi \big) \, dp
\\
-\int_{\Hcal_{0}} \phi(p)\, g \big(\langle \nu_{p}, F(\cdot,p) \rangle, \nn_{\Hcal_{0}} \big) \,dp + \int_{\Mbf} \phi(p)\langle\nu_{p}, U'(\cdot) \dive f(\cdot,p) \rangle \,dp \le 0.
\end{multline}
\end{lemma}

\begin{proof}[Proof of Theorem \ref{FVM-1}]
According to the inequality \eqref{POC.5}, we have for all convex entropy pairs
$(U,F)$,
$$
\dive  \langle \nu, F(\cdot) \rangle - \langle \nu, (\dive F)(\cdot) \rangle  + \langle \nu, U'(\cdot) (\dive f)(\cdot) \rangle \leq 0
$$
in the sense of distributions in $\Mbf$.
Since on the initial hypersurface $\Hcal_{0}$, the (trace of the) Young measure $\nu$ is the Dirac mass $\delta_{u_0}$ 
(because $u_0$ is a bounded function), from 
the theory in \cite{BL}
there exists a unique function $u \in L^\infty(\Mbf)$ such that the measure $\nu$ remains the 
Dirac mass $\delta_u$ for all Cauchy hypersurfaces $\Hcal_{t},$ $0\le t\le T$. Moreover, this implies that the approximations $u^h$ converge strongly to $u$
on compact sets at least. This concludes the proof.
\end{proof}

\begin{proof}[Proof of Lemma~\ref{POC-1}]
The proof consists of passing to the limit the inequality \eqref{FVM.31} and using the property \eqref{POC.1} of the Young measure.
First, note that the left-hand side of the inequality \eqref{FVM.31} converges immediately to the first integral term of \eqref{POC.5}, in view of \eqref{POC.1}. Next, take the last term of \eqref{FVM.31}. Using again \eqref{POC.1} and the fact that  $\phi_{\del^{0}K}^{n} - \phi(p) = \Ocal(\tau + h)$, we see that this term converges to the second integral term in \eqref{POC.5}.

Next, we will prove that the first term on the right-hand side of \eqref{FVM.31} converges to the last term in \eqref{POC.5}. Observe first that
$$
\tilde\mu^{n+1}_{K,\ekz} - \overline\mu^{n+1}_{K,\ekz} = \frac{1}{| e_{K^{n}}^{+}|} \int_{\Kn}\dive f(\ukn, p) dp.
$$
Therefore, we obtain
\begin{multline*}
-\sum_{n=0}^{\infty} \sum_{\sumkez} \frac{|\ekz|}{|\del^{0}\Kn|} |e_{K^{n}}^{+}| \phi_{\ekz}^{n} \big( V^{n+1}_{K}(\tilde\mu^{n+1}_{K,\ekz}) - V^{n+1}_{K}(\overline\mu^{n+1}_{K,\ekz}) \big)
\\
=\sum_{n=0}^{\infty} \sum_{\sumkez} \frac{|\ekz|}{|\del^{0}\Kn|}  \phi_{\ekz}^{n} \Big(
 \del_{\mu}V^{n+1}_{K}(\tilde\mu^{n+1}_{K,\ekz}) \int_{\Kn}\dive f(\ukn, p) dp) + | e_{K^{n}}^{+}|\Ocal(\tau^{2}) \Big)
\\
=\sum_{n=0}^{\infty} \sum_{\sumkez} \frac{|\ekz|}{|\del^{0}\Kn|}  \phi_{\ekz}^{n} \Big(
 \big( \del_{\mu}V^{n+1}_{K}(\tilde\mu^{n+1}_{K,\ekz}) - U'(\uknn) \big)\int_{\Kn}\dive f(\ukn, p) dp) 
\\
+ | e_{K^{n}}^{+}| \Ocal(\tau^{2}) + U'(\uknn) \int_{\Kn}\dive f(\ukn, p) dp) \Big).
\end{multline*}
Now, note that from the expression of $V$, \eqref{FVM.100},
\begin{align*}
\del_{\mu}V^{n+1}_{K}(\tilde\mu^{n+1}_{K,\ekz}) - U'(\uknn) &= U' ( (\muknn)^{-1}(\tilde\mu^{n+1}_{K,\ekz})) - U'(\uknn)
\\
&\le \sup U'' \sup_{n,\Kn}\Lip(\mukn)^{-1} | \tilde\mu^{n+1}_{K,\ekz} - \muknn|,
\end{align*}
and so, using the $L^{\infty}$ bound \eqref{FVM.12} and the growth condition \eqref{FVM.2}, we find
\be
\aligned 
& \sum_{n=0}^{\infty} \sum_{\sumkez} \frac{|\ekz|}{|\del^{0}\Kn|}  \phi_{\ekz}^{n} \big( \del_{\mu}V^{n+1}_{K}(\tilde\mu^{n+1}_{K,\ekz}) - U'(\uknn) \big)\int_{\Kn}\dive f(\ukn, p) \, dp
\\
& \lesssim \sum_{n=0}^{\infty} \sum_{\sumkez}\frac{|\ekz|}{|\del^{0}\Kn|}  \phi_{\ekz}^{n} |\Kn| | \tilde\mu^{n+1}_{K,\ekz} - \muknn|.
\endaligned
\ee
Applying Cauchy-Schwartz inequality and the entropy dissipation estimate \eqref{FVM.26}, we find that this term tends to zero with $h$. We are left with the term
$$
\sum_{n=0}^{\infty} \sum_{\Kn\in\Kcal^{n}} \phi_{\del^{0}K}^{n} U'(\uknn) \int_{\Kn}\dive f(\ukn, p) dp),
$$
which is easily seen to be of the form
$$
\sum_{n=0}^{\infty} \sum_{\Kn\in\Kcal^{n}} U'(\unk)\int_{\Kn} \phi(p) \dive f(\ukn,p) dp + \Ocal(h) \to  \int_{\Mbf}\langle\nu_{p}, U'(\cdot) \dive f(\cdot,p) \rangle \,dp.
$$

It remains to check that the remaining terms in \eqref{FVM.31} tend to zero with $h$. Namely, the second term on the right-hand side can be written as
\be
\aligned 
& \sum_{n=0}^{\infty} \sum_{\sumkez} \frac{|\ekz|}{|\del^{0}\Kn|} |e_{K^{n}}^{+}| (\phi_{\del^{0}K}^{n} - \phi_{\ekz}^{n}) V^{n+1}(\overline\mu^{n+1}_{K,\ekz})
\\
& = \sum_{n=0}^{\infty} \sum_{\sumkez} \frac{|\ekz|}{|\del^{0}\Kn|} |e_{K^{n}}^{+}| (\phi_{\del^{0}K}^{n} - \phi_{\ekz}^{n})\big( V^{n+1}(\overline\mu^{n+1}_{K,\ekz}) - V^{n+1}(\muknn(\uknn)) \big) 
\\
& = o(1),
\endaligned
\ee
by Cauchy-Schwartz inequality and the entropy dissipation estimate \eqref{FVM.26}. Next, the second term on the right-hand side of \eqref{FVM.31} satisfies
\be
\aligned 
& \sum_{n=0}^{\infty} \sum_{\sumkez} \int_{\ekz} (\phi_{\ekz}^{n} - \phi(p)) 
F_{\ekz}(\ukn,p) dp
\\
& = \sum_{n=0}^{\infty} \sum_{\sumkez} \int_{\ekz} (\phi_{\ekz}^{n} - \phi(p)) 
\big(F_{\ekz}(\ukn,p) - \dashint_{\ekz}F_{\ekz}(\ukn,q) dq \big) \, dp.
\endaligned 
\ee
In view of the regularity of $\phi$ and $F$, this term is bounded by
$$
\sum_{n=0}^{\infty} \sum_{\Kn\in\Kcal^{n}} |\del^{0}K| \Ocal(\tau_{\Kn} + h)^{2}.
$$
Using the CFL condition \eqref{FVM.7}, and the property \eqref{FVM.A}, we can further bound this term by
$$
\sum_{n=0}^{\infty} \sum_{\Kn\in\Kcal^{n}} | e_{K^{n}}^{+}| \Ocal(\tau_{\Kn})\{\Ocal(\tau_{\Kn}+h) + \Ocal(h^{2}/\tau_{\Kn})\} = o(1).
$$
Only the term 
\begin{multline}
\label{POC.10}
A^{h}(\phi) :=
\sum_{n=0}^{\infty} \sum_{\Kn\in\Kcal^{n}}\Bigg(
 \int_{e_{K^{n}}^{-}} (\phi_{\del^{0}K}^{n-1} - \phi(p))   \,g( F(\ukn,p) , \nn_{K^{n},\ekm}) dp 
\\
+ \int_{e_{K^{n}}^{+}} (\phi_{\del^{0}K}^{n} - \phi(p))  \,g( F(\ukn,p) , \nn_{K^{n},\ekp}) dp \Bigg) 
\end{multline}
remains to be bounded.
Here, we will use that our triangulation is admissible, in the sense of Definition \ref{FVM-0-2}.
First of all, observe that using Lemma \ref{PC-4-1}, one may replace $\phi_{\del^{0}K}^{n}$ by $\phi( p^{0}_{\Kn})$, and $\phi_{\del^{0}K}^{n-1}$ by $\phi( p^{0}_{K^{n-1}})$, where $ p^{0}_{K^{j}}$ denotes the center of $\del^{0}K^{j}$, with an error term of the form $C (\tau +h) \|\phi\|_{C^2} \| F \|_{L^{\infty}}$. Next, we replace 
$$
\phi(p)   \,g( F(\ukn,p) , \nn_{K^{n},\ekp})
$$
with 
$$
\phi( p^{+}_{K})  \,g( F(\ukn,\pkp) , \nn_{K^{n},\ekp}(\pkp)),
$$
and similarly for $\ekm$. Using the property \eqref{FVM.C}, the resulting error term is seen to be of the form $C(\tau+h) \|\phi\|_{C^2} \| F \|_{C^2}$. We then have
\begin{multline*}
A^{h}(\phi) \lesssim\sum_{n=0}^{\infty} \sum_{\Kn\in\Kcal^{n}} \Bigg(
 |\ekp| (\phi( p^{0}_{\Kn}) - \phi( p^{+}_{\Kn}) ) \, g \big( F(\ukn, p^{+}_{\Kn}), \nn_{K,\ekp}( p^{+}_{\Kn}) \big)
\\
+ |\ekm| (\phi( p^{0}_{K^{n-1}}) - \phi( p^{+}_{K^{n-1}}) ) \, g \big( F(\ukn, p^{+}_{K^{n-1}}) , \nn_{K,\ekm}( p^{+}_{K^{n-1}}) \big) \Bigg)
\\
 + (\tau+h) \|\phi\|_{C^2} \big( \| F \|_{L^{\infty}} + \| F \|_{C^2} \big).
\end{multline*}
Now, performing a Taylor expansion of $\phi$ and using the definition of $\wbf_K$ 
(recall that $\wbf_{K}$ is the future-oriented unit vector at $ p^{+}_{K}$ tangent to the geodesic connecting $ p^{+}_{K}$ and $ p^{0}_{K}$) we find, for instance,
$$
\phi( p^{0}_{K}) - \phi( p^{+}_{K}) \le h \, g(\wbf_{K}, \nabla \phi ( p^{+}_{K})) + \Ocal(h^{2}).
$$
Therefore, by the definition of $\Ecal(K)$ and 
using \eqref{FVM.E}, we may express this using the Cartesian deviation of the triangulation,
$$
\aligned
A^{h}(\phi) &\lesssim \sum_{n=0}^{\infty} \sum_{\Kn\in\Kcal^{n}} h \Bigg(
  |\ekp| \, g(\wbf_{\Kn}, \nabla \phi ( p^{+}_{K}) ) \, g( F(\ukn,  p^{+}_{K}), \nn_{\Kn,\ekp}( p^{+}_{\Kn}) ) 
\\
&\qquad+  |\ekm| \, g(\wbf_{K^{n-1}}, \nabla \phi ( p^{+}_{K^{n-1}}) ) \, g( F(\ukn,  p^{+}_{K^{n-1}}), \nn_{\Kn,\ekm}( p^{+}_{K^{n-1}}) ) \Bigg) 
\\
&\qquad+(\tau+h) \|\phi\|_{C^2}(\| F \|_{L^{\infty}} + \| F \|_{C^2} )
\\
&\le h \sum_{\Kn\in\Tcal^{h}} \big( |K| \Ecal(K) - |K^{-}| \Ecal(K^{-}) \big) (\nabla\phi, F(u^{h}) ) 
\\
&\qquad+ (\tau +h) \|\phi\|_{C^2}(\| F \|_{L^{\infty}} + \| F \|_{C^{2}} )
\\
&\lesssim \eta(h) \| \phi \|_{C^{1}} \| F \|_{L^{\infty}} +  (\tau +h) \|\phi\|_{C^{2}}(\| F \|_{L^{\infty}} + \| F \|_{C^{2}} ),
\endaligned
$$
which tends to zero since $\eta(h) \to 0$. This completes the proof of Lemma~\ref{POC-1}.
\end{proof}


\section{Examples and remarks}
\label{EX}

\subsection{A geometric condition}

The following proposition gives a condition of a geometric nature which is sufficient for a triangulation to be admissible in the sense of \eqref{FVM.E}. Recall that $ p^{\pm}_{K}$ denotes the center of mass of $\ekpm$, and that the vector $\wbf$ denotes the tangent at $\ekp$ to the geodesic from $\ekp$ to the center of $\del^{0}K$.
\begin{proposition}
\label{EX-1}
Let $\Tcal^{h}$ be a triangulation. Suppose that, for each element $K$, the scaled exterior normals $|\ekp|\nn_{K,\ekp}$ and $|\ekm| \nn_{K,\ekm}$, and the vectors $\wbf_{K}$ and $\wbf_{K^{-}}$ approach in the limit, in the following weak sense: for every smooth vector field $X$, one has
\be
\label{EX.3}
\big| g( |\ekp| \nn_{K,\ekp}, X) - g( |\ekm| \nn_{K,\ekm}, X)\big|  \lesssim \frac{\eta(h)}{h} \,|K| \| X\|_{L^{\infty}(K)},
\ee
\be
\label{EX.4}
\big| g( \wbf_{K}, X) - g( \wbf_{K^{-}}, X)\big|  \lesssim \frac{\eta(h)}{h} \,\tau_{K} \| X\|_{L^{\infty}(K)},
\ee
where the expressions are evaluated at the centers of mass of $\ekp$ and $\ekm$, and $\eta(h)$ is such that $ \eta(h) \to 0$. Then, $\Tcal^{h}$ is an admissible triangulation in the sense of \eqref{FVM.E}.
\end{proposition}
\begin{proof}
Let $\Phi$ be a smooth vector field, and let $\Psi$ be smooth on each element $K$. We have
$$
\aligned
&\sum_{K\in\Tcal^{h}} \big( |K| \Ecal(K) - |K^{-}| \Ecal(K^{-}) \big) (\Phi, \Psi _{/{K}}) 
\\
& = \sum_{K\in\Tcal^{h}} \big( |\ekp| \wbf_{K}\otimes \nn_{K,\ekp} - |\ekm| \wbf_{K^{-}}\otimes \nn_{K,\ekm} \big) (\Phi, \Psi _{/{K}}) 
\\
& = \sum_{K\in\Tcal^{h}} g(\wbf_{K}, \Phi) \, g(|\ekp|\nn_{K,\ekp}, \Psi_{/K}) -  g(\wbf_{K^{-}}, \Phi) \, g(|\ekm|\nn_{K,\ekm}, \Psi_{/K})
\\
& = \sum_{K\in\Tcal^{h}} g(\wbf_{K^{-}}, \Phi)\big( g(|\ekp|\nn_{K,\ekp}, \Psi_{/K}) - g(|\ekm|\nn_{K,\ekm}, \Psi_{/K}) \big) 
\\
&\qquad + \big( g(\wbf_{K}, \Phi) -g(\wbf_{K^{-}}, \Phi)\big) \, g(|\ekp|\nn_{K,\ekp}, \Psi_{/K})
\\
& \lesssim \eta(h)/h \sum_{K\in\Tcal^{h}} |K|  \| \Phi \|_{L^{\infty}} \| \Psi_{/K}\|_{L^{\infty}}
\lesssim \eta(h)/h \| \Phi \|_{L^{\infty}} \sup_{K} \| \Psi_{/K}\|_{L^{\infty}}.
\endaligned
$$
In view of \eqref{FVM.E}, this shows that $\Tcal^{h}$ is an admissible triangulation.
\end{proof}

For instance, one can easily check that if a triangulation is subordinate to a given foliation (in the sense that each connected component of the set of all outgoing faces 
$$\{ \ekp : K \in\Tcal^{h} \}$$
is contained in a certain Cauchy surface), and if, moreover, each lateral face $\ekz$ is everywhere tangent to a given, fixed, smooth time-like vector field, then the hypotheses of Proposition~\ref{EX-1} hold. However, our conditions \eqref{FVM.E} or \eqref{EX.3}, \eqref{EX.4} allow for more general triangulations, which need not satisfy such regularity assumptions.


\subsection{Lax--Friedrichs-type flux-functions}

Our general space-time setting requires some care when defining particular numerical fluxes. In particular, the absence of a canonically singled out conservative variable implies that one needs to use the values of $\mu^{f}_{K,e}(u)$, for both $K$ and its neighboring element $K_{\ekz}$, to compute the flux along an interface.
For instance, consider the following natural generalization of the Lax--Friedrichs-type flux,
\be
\label{EX.2}
\qq (u,v) = \frac{1}{2} \big(\mu^{f}_{K,\ekz}(u) + \mu^{f}_{K,\ekz}(v) \big) +
\frac{D_{K,\ekz}}{2} \big( \mu^{f}_{K,\ekp}(u) - \mu^{f}_{K_{\ekz},e^{+}_{K_{\ekz}}}(v)\big),
\ee
where the constants $D_{K,\ekz}$ satisfy
$$
D_{K,\ekz} \ge \frac{|\ekp|}{|\del^{0}K|}.
$$
This numerical flux is conservative, and it is monotone, as may be checked using the CFL condition \eqref{FVM.7}. The consistency property \eqref{FVM.3.1} is seen to be valid once we remark that this numerical flux can be written as
$$
\qq(u,v,\mu_{1}(u), \mu_{2}(v)) =  \frac{1}{2} \big(\mu^{f}_{K,\ekz}(u) + \mu^{f}_{K,\ekz}(v) \big) +
\frac{D_{K,\ekz}}{2} \big( \mu_{1}(u) - \mu_{2}(v)\big),
$$
and that the consistency property reads, in fact,
$$
\qq(u,u, \mu(u), \mu(u)) = \mu^{f}_{K,\ekz}(u).
$$


\section*{Acknowledgements}

The authors were partially supported by the Agence Nationale de la Recherche (A.N.R.)
through the grant 06-2-134423
entitled {\sl ``Mathematical Methods in General Relativity''} (MATH-GR), and by the Centre National de la Recherche Scientifique (CNRS). The first author
(P.A.) was also supported by the {\sl FCT--Funda\c c\~ ao para a Ci\^encia e Tecnologia} 
(Portuguese Foundation for Science and Technology) through the grant SFRH/BD/17271/2004. 

 
\newcommand{\auth}{\textsc}


\end{document}